\documentclass[a4paper,12pt]{amsart}

\usepackage{latexsym, amsfonts, amsmath, amsthm, amssymb, verbatim}
\newcommand{\AC}{{AC}}
\newcommand{\BV}{{BV}}

\newcommand{\mR}{\mathbb{R}}
\newcommand{\mC}{\mathbb{C}}

\newcommand{\Real}{{\rm Re}}
\newcommand{\Imag}{{\rm Im}}

\newcommand{\lam}{\boldsymbol{\lambda}}
\newcommand{\xvec}{\boldsymbol{x}}
\newcommand{\yvec}{\boldsymbol{y}}

\newcommand{\norm}[1]{\left\lVert#1\right\rVert}

\newcommand{\normbv}[1]{\left\lVert#1\right\rVert_{\BV(\sigma)}}

\newcommand{\newnorm}[1]{|\!|\!|\, #1\, |\!|\!|}

\renewcommand{\mod}[1]{\left \vert #1 \right \vert}
\newcommand\st{\thinspace : \thinspace}
\def\sbrack#1{\Bigl \{ #1 \Bigr \} }

\def\snorm#1{\Bigl \Vert #1 \Bigr \Vert}

\def\sparen(#1){\Bigl ( #1 \Bigr )}
\def\ssparen(#1){ (#1) }

\newcommand{\Proj}{{\rm Proj}}

\DeclareMathOperator{\var}{var}
\DeclareMathOperator{\cvar}{\rm cvar}
\DeclareMathOperator{\vf}{vf}

\newtheorem{thm}{Theorem}[section]
\newtheorem{cor}[thm]{Corollary}
\newtheorem{lem}[thm]{Lemma}
\newtheorem{prop}[thm]{Proposition}

\newtheorem{ex}[thm]{Example}
\newtheorem{rem}[thm]{Remark}

\newcommand{\booktitle}[1]{\textit{#1}}
\newcommand{\journalname}[1]{\textrm{#1}}

\title{Compact $\AC(\sigma)$ operators}
\author{Brenden Ashton and Ian Doust}
\address{School of Mathematics and Statistics\\
University of New South Wales\\
UNSW Sydney 2052 Australia}%
\email{i.doust@unsw.edu.au}%
\date{}
\begin{document}
\maketitle

%%%%%%%%%%%%%%%%%%%%%%%%%%%%%%%%%%%%%%%%%%%%%%%%%%%%%%%%%%%%%%%%%%%%%%%
%
%   The abstract
%
%%%%%%%%%%%%%%%%%%%%%%%%%%%%%%%%%%%%%%%%%%%%%%%%%%%%%%%%%%%%%%%%%%%%%%%

\begin{abstract}
All compact $\AC(\sigma)$ operators have a representation analogous
to that for compact normal operators. As a partial converse we
obtain conditions which allow one to construct a large number of
such operators. Using the results in the paper, we answer a number
of questions about the decomposition of a compact $\AC(\sigma)$ into
real and imaginary parts.
\end{abstract}

%%%%%%%%%%%%%%%%%%%%%%%%%%%%%%%%%%%%%%%%%%%%%%%%%%%%%%%%%%%%%%%%%%%%%%%
%
%   Introduction
%
%%%%%%%%%%%%%%%%%%%%%%%%%%%%%%%%%%%%%%%%%%%%%%%%%%%%%%%%%%%%%%%%%%%%%%%

\section{Introduction}
The class of well-bounded operators was introduced to provide a
theory which would allow many of the results which apply to
self-adjoint operators to be extended to the Banach space setting.
Since many operators which are self-adjoint on $L^2$ have only
conditionally convergent spectral expansions on the other $L^p$
spaces, the theory needed to allow more general types of
representation theorems than those available in the theory of
spectral operators.

The issue of the conditional convergence of spectral expansions
arises most explicitly when considering compact well-bounded
operators. In \cite{CD} it was shown that every compact well-bounded
operator $T$ has an expansion as
  \begin{equation}\label{main-sum}
   T = \sum_{j=1}^\infty \mu_j P_j
  \end{equation}
where $P_j$ is the Riesz projection onto the eigenspace
corresponding to the eigenvalue $\lambda_j$, and where the terms are
ordered so that $|\mu_j|$ is decreasing. Indeed, necessary and
sufficient conditions were given which ensure that any sum of the
form appearing in (\ref{main-sum}) is a compact well-bounded
operator.

Even in the earliest papers (see, for example, \cite{R2}) the fact
that the spectrum of a well-bounded operator is necessarily real was
seen as an undesirable restriction, and various attempts at
addressing this have appeared. In \cite{BG} Berkson and Gillespie
introduced the concept of an $\AC$ operator, which is one which can
be written in the form $T = A + i B$ where $A$ and $B$ are commuting
well-bounded operators. Doust and Walden \cite{DW} showed that, as
long as one takes the eigenvalues in an appropriate order, every
compact $\AC$ operator has a representation in the form given in
(\ref{main-sum}).
%The function theory associated to these
%operators was restricted to functions defined on a rectangle
%containing the spectrum.

The theory of $\AC$ operators had certain drawbacks however (see
\cite{BDG}) and a smaller class of operators, known as $\AC(\sigma)$
operators, was introduced in \cite{AD1}. Whereas the functional
calculus associated to theory of $\AC$ operators was restricted to
functions defined on a rectangle in the plane, the theory of
$\AC(\sigma)$ operators is based on functions of bounded variation
which are defined on an arbitrary nonempty compact subset $\sigma
\subseteq \mathbb{C}$.

The results of \cite{DW} clearly also apply to compact $\AC(\sigma)$ operators, but even in the case of $\AC$
operators, there has not been a characterization of the compact operators in terms of properties of the
eigenvalues and the corresponding projections. One of the main applications of the characterization result in
\cite{CD} has been to enable the construction of well-bounded operators with specific properties. See, for
example, \cite{CD,DG,DL}. The main aim of this paper, then, is to prove Theorem~\ref{construct-AC} which gives
sufficient conditions to ensure that an operator of the form (\ref{main-sum}) is a compact $\AC(\sigma)$
operator. To prove this theorem one has to show that (under the hypotheses of the theorem) one may sensibly
define $f(T)$ for $f \in \AC(\sigma)$, with a norm bound $\norm{f(T)} \le K \normbv{f}$.  A significant
challenge in working with $\AC(\sigma)$ operators is begin able to calculate $\normbv{f}$ for $f \in
\AC(\sigma)$. In Section~\ref{Norm-est} we shall show that for certain sets $\sigma$, the $\AC(\sigma)$ norm is
equivalent to a norm which is much easier to calculate. Although we will not need the full force of this result
to prove Theorem~\ref{construct-AC}, we feel that this equivalence is of independent interest.

In Section~\ref{S-ex} we show that there are compact $\AC(\sigma)$ operators which are not of the form
constructed in Theorem~\ref{construct-AC}. The final section includes a discussion of the properties of the
splitting of an $\AC(\sigma)$ operator $T$ into real and imaginary parts $T = A+iB$. There are many open
questions regarding these splittings. In the case that $T$ is compact however, it is possible to resolve these
questions. In obtaining these results we prove a new result about rearrangements of the sum representation of a
compact well-bounded operator (Corollary~\ref{same-sum}) which may also be of independent interest.

\section{Preliminaries}\label{S-preliminaries}

Throughout this paper let $\sigma \subset \mC$ be compact and
non-empty. For a Banach space $X$ we shall denote the bounded linear
operators on $X$ by $B(X)$ and the bounded linear projections on $X$
by $\Proj(X)$. Summations over empty sets of indices should always
be interpreted as having value zero.

The Banach algebra of functions of bounded variation on $\sigma$,
denoted $\BV(\sigma)$, was defined in \cite{AD1}. The norm in this
space is given by an expression of the form
  \begin{equation}\label{bvnorm}
   \normbv{f} = \sup_\gamma \cvar(f,\gamma) \, \rho(\gamma).
  \end{equation}
In (\ref{bvnorm}) the supremum is taken over all piecewise linear
curves $\gamma:[0,1] \to \mathbb{C}$ in the place. The term
$\cvar(f,\gamma)$ measures the variation of $f$ as one travels along
the curve $\gamma$, and $\rho(\gamma)^{-1}$ measures how `sinuous'
the curve $\gamma$ is. More precisely, the variation factor
$\vf(\gamma) \equiv \rho(\gamma)^{-1}$ is defined as the maximum
number of entry points of the curve $\gamma$ on any line in the
plane. (Heuristically one should think of this as the maximum number
of times any line intersects $\gamma$.) We refer the reader to
\cite{AD1} for the full definitions. The affine invariance of the
$\BV$ norm will be used repeatedly (with little comment) to pass
between estimates for functions on $\mR$ to those for functions
defined on other lines.

The closure of the polynomials in $z$ and $\overline{z}$ is the
subalgebra $\AC(\sigma)$ of absolutely continuous functions on
$\sigma$. An operator $T \in B(X)$ which admits an $\AC(\sigma)$
functional calculus is said to be an $\AC(\sigma)$ operator. This
class includes all well-bounded operators, as well as every normal
operator on a Hilbert space $H$.

\begin{thm}\label{can-use-DW} Suppose that $T \in B(X)$ is a compact $\AC(\sigma)$ operator.
Then
\begin{enumerate}
\item\label{BG-1} $T$ is an $\AC$
operator (in the sense of Berkson and Gillespie \cite{BG});
\item\label{BG-2} there exist unique commuting well-bounded operators $A,B \in B(X)$
such that $T = A + iB$;
\item\label{BG-3} $A$ and $B$ are compact.
\end{enumerate}
\end{thm}

\begin{proof}
Statement (\ref{BG-1}) is Theorem~5.3 of \cite{AD1}. Statements
(\ref{BG-2}) and (\ref{BG-3}) therefore follow  from
\cite[Theorem~6.1]{DW}.
\end{proof}

%\begin{thm}\cite[Theorem 6.1]{AD2}
%If $T \in B(X)$ is an $\AC(\sigma)$ operator then $T$ is an $\AC$
%operator (in the sense of Berkson and Gillespie \cite{BG}), and
%hence there exist commuting well-bounded operators $A,B \in B(X)$
%such that $T = A + iB$.
%\end{thm}

For a complex number $\mu = x+i y$ with $x,y \in \mathbb{R}$, let
$\mod{\mu}_\infty = \max\{x,y\}$.  We shall now define an order
$\prec$ on $\mathbb{C}$ by setting $\mu_1 \prec \mu_2$ if
\begin{itemize}
  \item[(i)] $\mod{\mu_1}_\infty > \mod{\mu_2}_\infty$, or,
  \item[(ii)] if $\mod{\mu_1}_\infty
= \mod{\mu_2}_\infty = \alpha$ and $\mu_2$ lies on the that part of
the square $\mod{z}_\infty = \alpha$ between $-\alpha + i \alpha$
and $\mu_1$ going from in a clockwise direction.
\end{itemize}

%The following result appears in \cite{DW} under the stronger
%hypothesis that $T$ is compact. This hypothesis was however only
%used to prove that $\sigma(T)$ was a countable set with $0$ as its
%only limit point.

Theorem~\ref{can-use-DW} has as an immediate corollary that compact
$\AC(\sigma)$ operators have a spectral diagonalization analogous to
that for compact normal operators, but where the sum in the
representation might only converge conditionally.

\begin{cor}\label{DW-result}\cite[Theorem 4.5]{DW}
 Suppose that $T$ is a compact $\AC(\sigma)$ operator
with spectrum $\{0\} \cup \{\mu_j\}_{j=1}^\infty$ and that
$\{\mu_j\}$ is ordered by $\prec$.  Then there exists a uniformly
bounded sequence of disjoint projections $P_j \in B(X)$ such that
  \[T = \sum_{j=1}^\infty \mu_j P_j,\]
where the sum converges in the norm topology of $B(X)$.
\end{cor}

This includes, for example, the fact that the range of the Riesz
projection associated with a nonzero eigenvalue $\mu$ is exactly the
corresponding eigenspace. We refer the reader to \cite{DW} for a
fuller discussion of properties of compact $\AC$ operators.

\section{Approximation in $\AC(\sigma)$}\label{S-approx}

An important step in proving Corollary~\ref{DW-result} is that the identity function $\lam(z) = z$ can be
approximated in $\BV$ norm by functions whose support intersects $\sigma(T)$ at only a finite number of points.
It follows from the results in \cite{AD1} and  \cite{AD2} that this is still true if one uses the $\BV(\sigma)$
norm. We provide here a more direct proof, the results of which we will use later.

Let $\sigma$ be a nonempty compact subset of $\mathbb{C}$. Given $r
> 0$ and $\epsilon > 0$, define $g_{r,\epsilon}: \sigma \to
\mathbb{C}$ by
  \[ g_{r,\epsilon}(z) = \begin{cases}
          0, & \hbox{if $|z|_\infty \le r$,} \\
          (|z|_\infty - r)/\epsilon,
             & \hbox{if $r < |z|_\infty \le r+\epsilon$,} \\
          1, & \hbox{if $|z|_\infty > r+\epsilon$.}
          \end{cases} \]

\begin{lem} For all $r > 0$ and $\epsilon > 0$, $g_{r,\epsilon} \in \AC(\sigma)$ and
$\normbv{g_{r,\epsilon}} \le 6$.
\end{lem}

\begin{proof} Let $\sigma_R = \mathrm{Re}(\sigma) = \{x \st x + iy \in
\sigma\}$ and let $\sigma_I = \mathrm{Im}(\sigma) = \{y \st x + iy
\in \sigma\}$. Define $u: \sigma_R \to \mathbb{C}$ and $v: \sigma_I
\to \mathbb{C}$ by  $u(x) = g_{r,\epsilon}(x)$ and $v(y) =
g_{r,\epsilon}(iy)$. Clearly $u \in \AC(\sigma_R)$ with
$\norm{u}_{BV(\sigma_R)} \le 3$ (and similarly for $v$). Let
$\hat{u}(x+iy) = u(x)$, and $\hat{v}(x+iy) = v(y)$. By
\cite[Proposition 4.4]{AD1}, $\hat{u}, \hat{v} \in \AC(\sigma)$ with
$\normbv{\hat{u}} \le 3$ and $\normbv{\hat{v}} \le 3$. Now it is
easy to check that $g_{r,\epsilon} = \hat{u} \land \hat{v}$, and
hence, by \cite[Proposition 2.10]{AD3}, $g_{r,\epsilon} \in
\AC(\sigma)$ and $\normbv{g_{r,\epsilon}} \le 6$.
\end{proof}

Let $\xvec(z) = \Real(z)$ and let $\yvec(z) = \Imag(z)$, so that
$\lam = \xvec+i\yvec$.

\begin{lem}\label{lam-approx} Suppose that $\{r_n\}$ and $\{\epsilon_n\}$ are
sequences of positive numbers which converge to $0$. For $n =
1,2,\dots$, let $g_n = g_{r_n,\epsilon_n}$, let $\xvec_n = g_n
\xvec$ and let $\yvec_n = g_n \yvec$. Then $\xvec_n \to \xvec$ and
$\yvec_n \to \yvec$ in $\AC(\sigma)$. Consequently $g_n \lam_n \to
\lam$ in $\AC(\sigma)$.
\end{lem}

\begin{proof} It suffices to show that $\xvec_n \to \xvec$.
Now $\xvec_n - \xvec = (1-g_n)\, \xvec = (1-g_n)\, {\tilde
x}_n$ where
  \[ {\tilde x}_n(z) = \begin{cases}
            \Real(z), & \hbox{if $|\Real(z)| \le r_n+\epsilon_n$,}\\
            0         & \hbox{if $|\Real(z)| > r_n+\epsilon_n$.}
            \end{cases}
            \]
Now $\normbv{{\tilde x}_n} \le 5 (r_n + \epsilon_n)$
 and so
  \[ \normbv{\xvec_n - \xvec} \le \normbv{1- g_n} \normbv{{\tilde x}_n}
     \le 30 (r_n + \epsilon_n) \to
              0\]
as $n \to \infty$.
\end{proof}

\begin{rem}{\rm
It is clear that one could replace $g_n$ in the above proof with many other families of `cut-off' functions.  In
the proof of \cite[Theorem~4.5]{DW}, for example, the cut-off functions are based on L-shaped regions rather
than squares.}
\end{rem}

\section{Norm estimates in $\AC(\sigma)$}\label{Norm-est}

In order to show that an operator $T$ admits an $\AC(\sigma)$
functional calculus, one often needs to find estimates for both
$\norm{f(T)}$ and $\normbv{f}$. This can be difficult, even for
quite simple functions.

If $\sigma$ lies inside the union of a finite number of lines through the origin, then we shall show that it is
possible to decompose $f \in \AC(\sigma)$ into a sum of simpler functions in a way that allows good estimation
of the norms. The main issue is the following. Suppose that $\mathrm{supp} f \subseteq \sigma_0 \subseteq
\sigma$ for some compact set $\sigma_0$. One always has that $\norm{f|\sigma_0}_{\BV(\sigma_0)} \le \normbv{f}$.
The challenge is to prove an estimate of the form $\normbv{f} \le C \norm{f|\sigma_0}_{\BV(\sigma_0)}$. Even if
$\sigma \subseteq \mR$ such an estimate need not exist, so any results need to rely on geometric properties of
$\sigma_0$ and $\sigma$.

We begin with a technical lemma. As in \cite{AD1}, given points
$z_1,z_2,\dots,z_k \in \mC$, let $\Pi(z_1,z_2,\dots,z_k)$ denote the
piecewise linear path with these points as vertices.
%In general it is not always easy to show that an operator $T$ has an
%$\AC(\sigma)$ functional calculus. This is in part due to the fact
%that it is difficult to estimate both $\norm{f(T)}$ and
%$\normbv{f}$, even for quite simple functions. In this section we
%will deal with the second of these issues. We begin with a technical
% lemma.

%In this section we shall give sufficient conditions on a sequence of
%scalars $\{\mu_j\}$ and a sequence of projections $\{P_j\}$ which
%will ensure that the sum $\sum_{j=1}^\infty \mu_j P_j$ converges and
%gives and $\AC(\sigma)$ operator.

\begin{lem}\label{rho-est}
Suppose that $k \ge 2$ and that $S = z_1,z_2,\dots,z_k$ is a list of
complex numbers such that no two consecutive numbers lie in the
complement of the real axis. Let
  \begin{align*}
  J_1 & = \{ j \st z_j,z_{j+1} \in \mathbb{R}\}, \\
  J_2 & = \{ j \st z_j \in \mathbb{R},\ z_{j+1} \not\in
               \mathbb{R}\}, \\
  J_3 & = \{ j \st z_j\not\in \mathbb{R},\ z_{j+1} \in
               \mathbb{R}\}
  \end{align*}
have cardinalities $k_1,k_2$ and $k_3$ respectively. Let $\gamma_S =
\Pi(z_1,z_2,\dots,z_k)$. Then
  \[ (k_2+k_3) \rho(\gamma_S) \le 2. \]
\end{lem}

\begin{proof} The conditions on $S$ imply that $|k_2-k_3| \le 1$.
The bound claimed obviously holds if $k_2 = k_3 = 0$, so we shall
assume that at least one of these values is nonzero.

Suppose first that $k_2 \le k_3$. If $j \in J_3$, then $z_{j+1}$ is
an entry point of $\gamma_S$ on the real axis. Thus $\rho(\gamma_S)
\le 1/k_3$ and so
  $ (k_2+k_3) \rho(\gamma_S) \le 2 k_3/k_3 = 2$.

If, on the other hand, $k_2 = k_3+1$, then the smallest element of
$J_2$ is less than the smallest element of $J_3$. Thus, in addition
to the entry points associated with the elements of $J_3$ (as in the
previous paragraph), $\gamma_S$ must have an earlier entry point on
the real line. Thus $\rho(\gamma_S) \le 1/(k+1)$, and so
  $ (k_2 + k+3) \rho(S) \le (2k_3+1)/(k_3+1) \le 2$.
\end{proof}

In general, if $\sigma_0 \subseteq \sigma$, then $\norm{f|\sigma_0}_{\BV(\sigma_0)} \le \normbv{f}$, but no
reverse inequality is available, even if $\mathrm{supp} f \subseteq \sigma_0$. Such an inequality does hold
however if $\sigma_0$ is a line inside $\sigma$.

Suppose then that $\sigma$ is a compact subset of $\mathbb{C}$ and that $\sigma_0 = \sigma \cap \mathbb{R} \ne
\emptyset$.

\begin{lem}\label{star-sigma-estimate} Suppose that $f \in \BV(\sigma)$ and that $\mathrm{supp} f
\subseteq \sigma_0$. Then
  \[ \normbv{f} \le 3\sbrack{\norm{f}_\infty + \var(f,\sigma_0)}
     = 3 \norm{f|\sigma_0}_{\BV(\sigma_0)}. \]
\end{lem}

\begin{proof}
For an ordered finite subset $S = \{z_1,\dots,z_k\}$ of $\sigma$
(allowing repetitions), let $\gamma_S = \Pi(z_1,\dots,z_k)$. Lemma
3.5 of \cite{AD1} shows that
  \[ \var(f,\sigma) = \sup_S \cvar(f,\gamma_S)\rho(\gamma_S) \]
where the supremum is taken over all such finite subsets. Indeed, by
adding extra points as necessary, one sees that
  \[ \var(f,\sigma) = \sup_S \sparen({\rho(\gamma_S)
       \sum_{j=1}^{k-1} |f(z_{j})-f(z_{j+1})|  }). \]
Fix such a subset $S$, and let $v(f,S) = \sum_{j=1}^{k-1}
  |f(z_{j})-f(z_{j+1})|$. Clearly $v(f,S)$ is unchanged if we omit
any consecutive elements of $S$ which are both in $\sigma\setminus
\sigma_0$. Note that omitting points never decreases the value of
$\rho(\gamma_S)$, so we shall assume that no two consecutive
elements of $S$ are both in $\sigma\setminus\sigma_0$. We may also
assume that $S \cap \mR \ne \emptyset$ (or else $v(f,S) = 0$).

Partition $J = \{1,2,\dots,k-1\}$ into sets $J_1,J_2,J_3$ as in
Lemma~\ref{rho-est}. Clearly then
  \[ \sum_{j=1}^{k-1} |f(z_{j})-f(z_{j+1})|
    = \sum_{i=1}^3 \sum_{j \in J_i} |f(z_{j})-f(z_{j+1})|. \]
Let $S_0 = \{w_1,\dots,w_{k_0}\}$ be the sublist of $S$ containing
the elements that lie on the real axis, and let $\gamma_{S_0}$
denote the corresponding piecewise linear curve. As noted above,
$\rho(\gamma_S) \le \rho(\gamma_{S_0})$. Thus, using
\cite[Proposition 3.6]{AD1},
  \[ \rho(\gamma_S) \sum_{j \in J_1} |f(z_{j})-f(z_{j+1})|
    \le \rho(\gamma_{S_0})
      \sum_{j=1}^{k_0-1} |f(w_{j})-f(w_{j+1})|
      \le \var(f,\sigma_0). \]
On the other hand, if $j \in J_2 \cup J_3$, then $| f(z_j) -
f(z_{j+1})| \le \norm{f}_\infty$ and so, by Lemma~\ref{rho-est},
  \[ \rho(\gamma_S) \sum_{j \in J_2 \cup J_3} |f(z_{j})-f(z_{j+1})|
     \le (k_2+k_3) \norm{f}_\infty \rho(\gamma_S)
     \le 2 \norm{f}_\infty. \]
Thus
  \[ \var(f,\sigma) \le \var(f,\sigma_0) + 2 \norm{f}_\infty \]
and hence
  \[ \normbv{f} \le 3 \norm{f}_\infty + \var(f,\sigma_0). \]
\end{proof}

Note that the factor 3 in the above inequality is sharp. If $\sigma
= \{i,0,-i\}$ and $f = \chi_{\{0\}}$ then $\normbv{f} = 3$, and
$\var(f,\{0\}) = 0$.

For $0 \le \theta < 2\pi$, let $R_\theta$ denote the ray $\{r\cos \theta,r \sin \theta) \st r \ge 0\}$. We shall
say that $\sigma \subseteq \mC$ is a \emph{spoke set} if it is a subset of a finite union of such rays.

Suppose then that $\sigma$ is a nonempty compact spoke set with $\sigma \subseteq \cup_{n=1}^N R_{\theta_n}$.
(We shall assume that the angles $\theta_n$ are distinct.) For $n = 1,\dots,N$, let $\sigma_n = \sigma \cap
R_{\theta_n}$. Given $f \in \BV(\sigma)$, we shall define $f_0,f_1,\dots,f_n \in \BV(\sigma)$ by setting $f_0(z)
= f(0)$ and, for $1 \le n \le N$,
  \[
  f_n(z) = \begin{cases}
     f(z)-f(0), & \hbox{if $z \in R_{\theta_n}$},\\
     0,    & \hbox{otherwise.}
     \end{cases}
  \]
Then $f = \sum_{n=0}^N f_n$. Also, if $f \in \AC(\sigma)$, then a short limiting argument can be used to show
that each $f_n$ is also in $\AC(\sigma)$. Define the spoke norm
  \[ \newnorm{f}_{Sp} = |f(0)| + \sum_{n=1}^N \norm{f_n|\sigma_n}_{\BV(\sigma_n)}. \]
Since $\sigma_n$ is lies in a line, the affine invariance of these norms and \cite[Proposition 3.6]{AD1} means
that calculating $\norm{f_n|\sigma_n}_{\BV(\sigma_n)}$ is relatively easy, since this just requires an
estimation of the usual variation of $f$ along the line. Thus $\newnorm{f}_{Sp}$ is much easier to calculate
than $\normbv{f}$.

The following result shows that $\newnorm{\cdot}_{Sp}$ and $\normbv{\cdot}$ are equivalent. Although we do not
need the full strength of this result in the next section, it does provide a useful tool in working with sets of
this sort.

\begin{prop}\label{equiv-norms} Suppose that $\sigma$ is a nonempty compact spoke set. Then for all $f \in \BV(\sigma)$
  \[ \frac{1}{2N+1}\, \newnorm{f}_{Sp} \le \normbv{f} \le 3\, \newnorm{f}_{Sp}. \]
\end{prop}

\begin{proof} Suppose that $1 \le n \le N$. Then
  \begin{align*}
  \norm{f_n}_{\BV(\sigma_n)}
     &= \sup_{\sigma_n} | f - f_0 | + \var(f-f_0,\sigma_n)  \\
     &\le 2 \norm{f}_\infty + \var(f,\sigma_n) \\
     &\le 2 \normbv{f}.
  \end{align*}
The left hand inequality then follows from the triangle inequality. On the other hand,
Lemma~\ref{star-sigma-estimate} (and affine invariance) shows that
  \[ \normbv{f} \le |f(0)| + \sum_{n=1}^n \normbv{f_n}
                \le |f(0)| + 3 \sum_{n=1}^n \norm{f_n}_{\BV(\sigma_n)}
                \le 3 \newnorm{f}_{Sp}. \]
\end{proof}

Let $bv_0$ denote the Banach space of sequences of bounded variation and limit $0$. The following lemma is
elementary.

\begin{lem} Suppose that $\{Q_j\}_{j=0}^\infty$ is a uniformly
bounded increasing family of projections on $X$. That is, $Q_i Q_j = Q_j Q_i = Q_i$ whenever $i \le j$ and
$\sup_j \norm{Q_j} \le K$. Suppose that $\{\mu_j\}_{j=1}^\infty \in bv_0$. Then
\begin{enumerate}
\item
$\snorm{\sum_{j=n}^m \mu_j (Q_j - Q_{j-1})} \le K \sparen(|\mu_n| + |\mu_m|
           +\sum_{j=n}^{m-1} |\mu_j - \mu_{j+1}|)$.
\item
$\sum_{j=1}^\infty \mu_j (Q_j - Q_{j-1})$ converges in norm.
\end{enumerate}
\end{lem}

\section{Constructing compact $\AC(\sigma)$ operators}

In \cite{CD3}, Cheng and Doust showed that certain combinations of disjoint projections of the form $\sum
\lambda_j E_j$ always converge and define compact real $\AC(\sigma)$ operators. In this section we shall provide
some sufficient conditions for an operator of this form to be a compact $\AC(\sigma)$ operator.
Theorem~\ref{construct-AC} below will allow the construction of compact $\AC(\sigma)$ operators which are
neither scalar-type spectral, nor real $\AC(\sigma)$ operators, via a given conditional decomposition of the
Banach space $X$.

%For $0 \le \theta < 2\pi$, let $R_\theta$ denote the ray $\{r\cos \theta,r \sin \theta) \st r > 0\}$.
Suppose
that $N \ge 1$ and that $\theta_1,\dots,\theta_n$ are distinct angles. Given
\begin{itemize}
  \item scalars $\{\lambda_{n,m} \st n=1,\dots,N,\ m=1,2,\dots \} \subset
\mathbb{C}$, and
  \item projections $\{E_{n,m} \st n=1,\dots,N,\ m=1,2,\dots \} \subset
B(X)$
\end{itemize}
consider the following three conditions:
\begin{enumerate}
  \item[(H1)] For each $n = 1,\dots,N$,  $\{\lambda_{n,m}\}_{m=1}^{\infty} \subset R_{\theta_n}$.
  \item[(H2)] For each $n = 1,\dots,N$, $|\lambda_{n,1}| >
  |\lambda_{n,2}| > |\lambda_{n,3}| > \dots$, and $\lambda_{n,m} \to 0$ as $m \to \infty$.
  \item[(H3)] The operators $E_{n,m}$ are pairwise disjoint, finite rank
  projections and there exists a constant $K$ such that,
   for each $n = 1,\dots,N$, and $M = 1,2,\dots$,
     $\displaystyle \snorm{\sum_{m=1}^M E_{n,m}}
   \le K$.
\end{enumerate}
The set of indices $\mathbb{I} = \{(n,m)\,:\, n=1,\dots,N,\ m=1,2,\dots\}$ can be ordered by declaring that
$(n,m) \succ (s,t)$ if $|\lambda_{n,m}| < |\lambda_{s,t}|$, or if $|\lambda_{n,m}| = |\lambda_{s,t}|$ and
$\theta_n > \theta_s$.

Let $\sigma = \{0\} \cup \{\lambda_{n,m}\}_{(n,m) \in \mathbb{I}}$, so that $\sigma$ is a spoke set.

\begin{thm}\label{construct-AC}
Suppose that $\{\lambda_{n,m} \st n=1,\dots,N,\ m=1,2,\dots\} \subset \mathbb{C}$ and $\{E_{n,m} \st
n=1,\dots,N,\ m=1,2,\dots\} \subset B(X)$ satisfy (H1), (H2) and (H3). Then
  \[ T = \sum_{n,m} \lambda_{n,m} E_{n,m} \]
converges in operator norm (in the order $\succ$) to a compact
$\AC(\sigma)$ operator.
\end{thm}

\begin{proof} Define $\Psi: \AC(\sigma) \to B(X)$ by
  \begin{equation}\label{Psi-defn}
   \Psi(f) = f(0)I + \sum_{n,m} (f(\lambda_{n,m}) - f(0)) E_{n,m}.
  \end{equation}
The first thing to verify is that $\Psi$ is well-defined, that is,
that the sum on the right-hand side of (\ref{Psi-defn}) converges
for all $f \in \AC(\sigma)$.

Suppose then that $f \in \AC(\sigma)$. Let $\mu_{n,m} = f(\lambda_{n,m}) - f(0)$. Fix $\epsilon > 0$. As $f$ is
continuous at $0$, if $(n_0,m_0)$ is large enough, then
  \begin{equation}\label{est-1}
  |\mu_{n,m}| < \epsilon/4N, \qquad \hbox{for all $(n,m) \succeq (n_0,m_0)$.}
  \end{equation}
As $f$ is absolutely continuous, for every $n = 1,2,\dots, N$, the variation along $R_\theta$,
$\var(f|\sigma_n)$, is finite. Hence, if $m_0$ is large enough,
  \begin{equation}\label{est-2}
  \sum_{m=s}^{t}
|\mu_{n,m}-\mu_{n,m+1}| < \epsilon/2N, \qquad \hbox{whenever $m_0 \le s \le t$.}
  \end{equation}

As $N$ is finite we can choose $n_0 \in \{1,\dots,N\}$ and $m_0 \ge 1$ such that (\ref{est-1}) holds and such
that (\ref{est-2}) holds for all $n$ at once. Suppose then that
 $(n_1,m_1) \succ (n_0,m_0)$. For each $n$, let $I_n$ be the (possibly empty) set
 $I_n = \{ m \,:\, (n_0,m_0) \preceq (n,m) \prec (n_1,m_1)\}$. If $I_n \ne \emptyset$, let $s_n = \min I_n$
 and $t_n = \max I_n$.
 The difference in the partial sums from index $(n_0,m_0$ to index $(n_1,m_1)$ is therefore given by
 \[ \Delta =  \sum_{{n \atop I_n \ne \emptyset}} \sum_{m=s_n}^{t_n} \mu_{n,m} E_{n,m}. \]
 Thus, by the
lemma,
  \begin{align*}
  \norm{\Delta}
    &\le \sum_{I_n \ne \emptyset}\ \snorm{\sum_{m=s_n}^{t_n} \mu_{n,m}
E_{n,m}} \\
    &\le \sum_{I_n \ne \emptyset} K \sparen( |\mu_{n,s_n}| + |\mu_{n,t_n}|
        + \sum_{m=s_n}^{t_n-1} |\mu_{n,m}-\mu_{n,m+1}|)\\
    &< K \epsilon
  \end{align*}
by (\ref{est-1}) and  (\ref{est-2}).  It follows that the partial sums are
  Cauchy and hence the series converges.  Note that in particular,
  this implies that the sum defining $T = \Psi(\lam)$
  converges. Since each $E_{n,m}$ is finite rank, $T$ is
  compact.

It is clear that $\Psi$ is linear. For $1 \le n \le N$, let $\sigma_n = \sigma \cap R_{\theta_n} = \{0\} \cup
\{\lambda_{n,m}\}_{m=1}^\infty$ as in Section~\ref{Norm-est}, and define $f_0,f_1,\dots,f_N$ as before
Proposition~\ref{equiv-norms}. Note that, using the affine invariance of $\AC(\sigma)$ and Lemma~3.2 of
\cite{CD3},
  \[ \norm{\Psi(f_n)}
  \le K \norm{f_n|\sigma_n}_{\BV(\sigma_n)}, \qquad \hbox{for $1 \le n \le N$}. \]
Then, by Proposition~\ref{equiv-norms},
  \begin{align*}
   \norm{\Psi(f)} &\le \sum_{n=0}^N \norm{\Psi(f_n)} \\
       &\le |f(0)| + K \sum_{n=1}^N \norm{f_n|\sigma_n}_{\BV(\sigma_n)} \\
       & \le K \newnorm{f}_{Sp} \\
       &\le (2N+1)K \normbv{f}.
   \end{align*}
It is easy to verify that $\Psi(fg) = \Psi(f)\Psi(g)$ if $f$ and $g$
are constant on a disk around $0$. The continuity of $\Psi$ then
implies that $\Psi$ is multiplicative on $\AC(\sigma)$.

Finally, since $T = \Psi(\lam)$, it follows that $T$ has an
$\AC(\sigma)$ functional calculus, and hence that $T$ is compact
$\AC(\sigma)$ operator.
\end{proof}

\section{Examples}\label{S-ex}

As one might expect, Theorem~\ref{construct-AC} is far from giving a
characterization of compact $\AC(\sigma)$ operators. Here we shall
give some examples which show that there are many ways of producing
compact $\AC(\sigma)$ operators whose spectra do not lie in a finite
number of lines through the origin.

\begin{prop}\label{construct-T} Let $\sigma = \{0,\lambda_1,\lambda_2,\dots\}$ be a
countable set of complex number whose only limit point is $0$.
Define $T$ on $\AC(\sigma)$ by $Tf(z) = z\,f(z)$. Then $T$ is a
compact $\AC(\sigma)$ operator.
\end{prop}

\begin{proof}
That $T$ has the required functional calculus is an immediate
consequence of that fact that $\AC(\sigma)$ is a Banach algebra. Let
$r_n = \epsilon_n = \frac{1}{n}$ and define $\lam_n = g_n \lam$ as
in Lemma~\ref{lam-approx}. It follows that $T = \lim_{n \to \infty}
\lam_n(T)$. But $\lam_n(T)$ is a finite rank operator, and hence $T$
is compact.
\end{proof}

An important class of examples is given \cite[Example 3.9]{AD3}.

\begin{prop} Let $A$ be a closed operator on a Banach space $X$, and
suppose that for some $x \in \rho(A)$, the resolvent $(xI - A)^{-1}$
is compact and well-bounded. Then $(\mu I - A)^{-1}$ is a compact
$\AC(\sigma_\mu)$ operator for all $\mu \in \rho(A)$.
\end{prop}

\begin{proof} Let $R(\mu,A) = (\mu I - A)^{-1}$. The resolvent
identity clearly implies that if one resolvent is compact then every
resolvent is compact. If we fix $\mu \not\in \sigma(T)$, then
$R(\mu,A) = f(R(x,A))$ where $f(t) = t/(1+(\mu-x)t)$ is a M{\"
o}bius transformation. If $J$ is any compact interval containing
$\sigma(R(x,A))$ then $\rho(f(J)) = \frac{1}{2}$. Thus $R(\mu,A)$ is
an $AC(f(J))$ operator.
\end{proof}

It was shown in \cite{DG} that there exist compact $\AC$ operators
(in the sense of Berkson and Gillespie) for which the sum
(\ref{main-sum}) fails to converge if the eigenvalues are listed in
order of decreasing modulus. It is not clear however whether that
construction always produces an $\AC(\sigma)$ operator.
%In this
%section we shall adapt the construction to produce a compact
%$\AC(\sigma)$ operator with this same property.
%We shall now construct a compact $\AC(\sigma)$ operator for which
%the sum $\sum_{j=1}^n \lambda_j P_j$ fails to converge when the
%eigenvalues are listed in order of decreasing modulus.
%The plan of
%the construction is close to that in \cite{DG},
The following adaption of that construction does produce an
$\AC(\sigma)$ operator example. Certain aspects require more care
here however due to the nature of the $\BV(\sigma)$ norm.

\begin{ex}{\rm Let $\theta = \tan^{-1}(1/6)$. For $k = 1,2,\dots$, let
  \begin{align*}
   \lambda_{k,j} & = \frac{e^{i\, j\theta/k}}{k}, & j=0,1,\dots,k,\\
   \mu_{k,j} & = \frac{\lambda_{k,j}+\lambda_{k,j-1}}{2}, & j =
   1,2,\dots,k.
   \end{align*}
Thus $d_k = |\mu_{k,j}|$ is independent of $j$, and $d_k <
|\lambda_{k,j}| = \frac{1}{k}$ for all $k,j$.

Let $\sigma = \{\lambda_{k,j}\}_{k,j} \cup \{\mu_{k,j}\}_{k,j} \cup
\{0\}$. Then $\sigma$ is compact, and hence by
Proposition~\ref{construct-T}, the operator $T \in B(\AC(\sigma))$,
$Tf(z) = zf(z)$ is a compact $\AC(\sigma)$ operator. Thus $T =
\sum_{\lambda_j=1}^\infty \lambda_j P(\lambda_j)$ where
$\{\lambda_j\}$ is a listing of the nonzero elements of $\sigma$
according to the order $\prec$ defined in
Section~\ref{S-preliminaries}, and $P(\lambda_j)$ is the projection
$P(\lambda_j) f = \chi_{\{\lambda_j\}} f$.

For $r > 0$, let $S_r = \sum_{|\lambda_j|\ge r} \lambda_j
P(\lambda_j)$, so that $S_r$ is a partial sum of the above series
for $T$ when the terms are ordered according to modulus. We shall
show that the series does not converge in this order by showing that
this sequence of partial sums is not Cauchy.

Fix $k$. Then
  \begin{align*}
   S_{d_k} - S_{1/k} &= \sum_{j=0}^k \lambda_{k,j} P(\lambda_{k,j})
   \\
     & = \lambda_{k,0} \sum_{j=0}^k P(\lambda_{k,j})
            + \sum_{j=1}^k (\lambda_{k,j}-\lambda_{k,0})P(\lambda_{k,j})
  .\end{align*}
Now $\lambda_{k,0} = 1/k$. Elementary trigonometry ensures that for
all $j$,
  \[ |\lambda_{k,j} - \lambda_{k,0}|
   \le |\lambda_{k,k} - \lambda_{k,0}|
   \le \frac{1}{k} \tan \theta = \frac{1}{6k}. \]
As $\AC(\sigma)$ is a Banach algebra, $\norm{P(\lambda_{k,j})} =
\normbv{\chi_{\{\lambda_{k,j}\}}} \le 3$ for all $j$. Thus
  \begin{equation}\label{small-bits}
   \norm{\sum_{j=1}^k
  (\lambda_{k,j}-\lambda_{k,0})P(\lambda_{k,j})}
  \le \frac{1}{2}.
  \end{equation}
Now $\sum_{j=0}^k P(\lambda_{k,j})$ is the projection of
multiplication by the characteristic function of the set $\Lambda_k
= \{\lambda_{k,0},\dots,\lambda_{k,k}\}$ and so
  \[ \norm{\sum_{j=0}^k P(\lambda_{k,j})}
    = \normbv{\chi_{\Lambda_k}}. \]
Let $\gamma_k$ denote the piecewise linear curve in $\mathbb{C}$
joining the elements of $\Lambda_k$ in order. Note that $\gamma_k$
passes through each of the points $\mu_{k,j}$. Clearly any line in
the plane has at most two entry points on $\gamma_k$ and so
$\rho(\gamma_k) = 1/2$. Thus
  \[ \cvar(\chi_{\Lambda_k},\gamma_k) \rho(\gamma_k)
      = 2(k-1) \, \frac{1}{2} = k-1\]
and so
  \[ \normbv{\chi_{\Lambda_k}} = \norm{\chi_{\Lambda_k}}_\infty
     + \sup_\gamma \cvar(\chi_{\Lambda_k},\gamma) \rho(\gamma)
     \ge 1 + (k-1) = k. \]
Thus, using (\ref{small-bits}),
  \begin{align*}
   \norm{S_{d_k} - S_{1/k}}
   & \ge \frac{1}{k} \norm{\sum_{j=0}^k P(\lambda_{k,j})} -
          \norm{\sum_{j=1}^k
  (\lambda_{k,j}-\lambda_{k,0})P(\lambda_{k,j})} \\
  &\ge 1 - \frac{1}{2} = \frac{1}{2}.
  \end{align*}
It follows that the partial sum sequence is not Cauchy and hence the
infinite sum does not converge. }
\end{ex}

\section{Other properties}

As was noted in Section~\ref{S-preliminaries}, every $\AC(\sigma)$
operator $T$ admits a splitting into real and imaginary parts $T =
A+iB$, where $A$ and $B$ are commuting well-bounded operators. On
nonreflexive spaces this splitting might not unique \cite[Example
4.5]{AD3}. Even on nonreflexive spaces however, one does get a
unique splitting when $T$ is compact. As is shown in the proof of
\cite[Theorem 6.1]{BDG}, this is because in this case the real and
imaginary parts are determined by the family of Riesz projections
associated with the nonzero eigenvalues of $T$. That is, if $T =
\sum \mu_j P_j$, then
  \[ A = \sum_i x_i \Bigl( \sum_{\Real(\mu_j) = x_i} P_j \Bigr) \]
where $\{x_i\}$ is the set of nonzero real parts of eigenvalues of
$T$, ordered so that $|x_1| \ge |x_2| \ge \dots$.

Given an $\AC(\sigma)$ operator $T$ and $\omega = \alpha + i \beta
\in \mC$, the operator $\omega T$ is an $\AC(\omega \sigma)$
operator. A longstanding open question is whether a splitting of
$\alpha T$ must be given by
  \begin{equation}\label{omega-T}
   \omega T = (\alpha A - \beta B) + i (\alpha B + \beta A).
  \end{equation}
Of course, if every real linear combination of $A$ and $B$ is again
well-bounded, then this must be a splitting. It is well-known
however that the sum of two commuting well-bounded operators $A$ and
$B$ need not be well-bounded (see \cite{G}). %, but no example is
%known where $A+iB$ is $\AC(\sigma)$ operator (rather than an
%$\AC$~operator.

If $T$ is compact, then $\omega T$ is also obviously compact, and
hence has a unique splitting as
  $ \alpha T = U + iV$.
The main issue in showing that (\ref{omega-T}) holds is in
rearranging the conditionally convergent sums that arise. The
following lemma shows that while rearrangements of the sum in
(\ref{main-sum}) may fail to converge, they cannot converge to a
different limit.

\begin{lem}\label{wb-rearrange}
Suppose that $\{c_j\}_{j=1}^\infty$ is a sequence of real numbers
whose only limit point is $0$. Suppose that $\{P_j\}$ is a sequence
of disjoint finite rank projections and that there is a constant $K$
such that
\begin{itemize}
\item for all $t  > 0$, $\snorm{\sum_{c_j \ge t} P_j} \le K$,
\item for all $t  < 0$, $\snorm{\sum_{c_j \le t} P_j} \le K$.
\end{itemize}
Then $\sum_{j=1}^\infty c_j P_j$ is well-bounded if the sum
converges.
\end{lem}

\begin{proof} Without loss we may assume that no $c_j$ is zero.
Suppose that $U = \sum_{j=1}^\infty c_j P_j$ converges. Let $\sigma
= \{0\} \cup \{c_j\}$. Then $U$ is clearly a compact operator with
$\sigma \subseteq \sigma(U)$. If $\beta \in \sigma(U)\setminus
\sigma$ then it is an isolated eigenvalue with corresponding Riesz
projection $P_\beta$. But then $UP_\beta = \beta P_\beta = \sum c_j
P_j P_\beta = 0$ which is impossible. Thus $\sigma(U) = \sigma$. It
is easy to confirm that the Riesz projection corresponding to $c_j$
is $P_j$.

Let $\pi$ be a permutation of the positive integers so that
$|c_{\pi(j)}|$ is nonincreasing. It follows from
\cite[Theorem~3.3]{CD} that $V = \sum_{j=1}^\infty c_{\pi(j)}
P_{\pi(j)}$ converges to a well-bounded operator (with $\sigma(V) =
\sigma$).

Let $\AC_{c}(\sigma) = \{ f \in \AC(\sigma) \,:\, \hbox{$f$ is
constant on a neighbourhood of $0$}\}$. Then $\AC_{c}(\sigma)$ is
dense in $\AC(\sigma)$. Let $\mathcal{A}$ denote the algebra of
functions $f$ which are analytic on a neighbourhood of $\sigma$, and
for which the restriction of $f$  to $\sigma$ lies in
$\AC_{c}(\sigma)$. Note that every $f \in \AC_{c}(\sigma)$ has an
extension to a locally constant element of $\mathcal{A}$. Suppose
then that $f \in \mathcal{A}$. Write $\sigma = \sigma_1 \cup
\sigma_2$, where $\sigma_2$ is the component of the spectrum
containing $0$ on which $f$ is constant, and where its complement
$\sigma_1$ is finite. The Riesz functional calculus for $U$ and $V$
gives that
  \[ f(U) = \sum_{c_j \in \sigma_1} f(c_j) P_j
      + f(0)\Bigl( I - \sum_{c_j \in \sigma_1} P_j \Bigr)
      = f(V). \]
But $V$ is well-bounded and so $\norm{f(U)} = \norm{f(V)} \le K
\normbv{f}$ for some $K$. The density of $\AC_{c}(\sigma)$ now
implies that $U$ is well-bounded.
\end{proof}

Note in particular that in the above proof,  if $f_n \in
\AC_{c}(\sigma)$, and $f_n \to \lam$ in $\AC(\sigma)$, then $U =
\lim_n f_n(U) = \lim_n f_n(V) = V$. This proves the following
result.

\begin{cor}\label{same-sum}
Suppose that $T$ is a compact well-bounded operator with
sum representation $\sum_j \mu_j P_j$ with $|\mu_1| \ge |\mu_2| \ge
\dots$. Let $\pi$ be a permutation of the positive integers. Then
  \[ T = \sum_{j} \mu_{\pi(j)} P_{\pi(j)} \]
if the sum on the right-hand side converges.
\end{cor}

It might be noted that we have been unable to prove the
corresponding result for compact $\AC(\sigma)$ operators.

We return now to the question raised at the beginning of this
section.

\begin{thm}
Let $T$ be a compact $\AC(\sigma)$ operator with splitting $T = A +
i B$, and let $\omega = \alpha + i \beta \in \mC$. The unique
splitting of $\alpha T$ is
  \[
   \omega T = (\alpha A - \beta B) + i (\alpha B + \beta A).
  \]
\end{thm}

\begin{proof} Write $T = \sum_{\succ} \mu_j P_j$ via
Corollary~\ref{DW-result}. Let $\xvec  = \Real(\lam)$ and $\yvec  =
\Imag(\lam)$. The proof of Corollary~\ref{DW-result} (see
Section~\ref{S-approx}) shows that the sums $\xvec(T) = \sum_{\succ}
\Real(\mu_j) P_j$ and $\yvec(T) = \sum_{\succ} \Imag(\mu_j) P_j$
both converge. Thus
  \[ T = A+iB = \sum_{\succ} \Real(\mu_j)
     P_j + i  \sum_{\succ} \Imag(\mu_j) P_j\]
and so
  \begin{align}
  \omega T &= \alpha \sum_{\succ} \Real(\mu_j) P_j
     + i \beta \sum_{\succ} \Real(\mu_j) P_j
     + i \alpha \sum_{\succ} \Imag(\mu_j) P_j
     - \beta \sum_{\succ} \Imag(\mu_j) P_j \notag\\
  &= \sum_{\succ} \Real(\omega \mu_j) P_j
                + i \sum_{\succ} \Imag(\omega \mu_j) P_j. \label{omegaT-splitting}
  \end{align}
The $\AC(\sigma)$ functional calculus for $T$ now provides the
bounds on the norms of sums of the Riesz projections needed to so
that we may apply Lemma~\ref{wb-rearrange} and deduce that
$\sum_{\succ} \Real(\omega \mu_j) P_j$ and $\sum_{\succ}
\Imag(\omega \mu_j) P_j$ are well-bounded. Since these operators
clearly commute, Equation~(\ref{omegaT-splitting}) gives the unique
splitting of $\omega T$. But $\sum_{\succ} \Real(\omega \mu_j) P_j =
\alpha A - \beta B$ and $\sum_{\succ} \Imag(\omega \mu_j) P_j =
\alpha B + \beta A$ so the proof is complete.
\end{proof}

The known examples of $\AC$ operators which are not $\AC(\sigma)$
operators share the property that they can be written as $T = A+iB$
where $A$ and $B$ are commuting well-bounded operators whose sum is
not well-bounded. The previous theorem shows that, at least for
compact operators, the well-boundedness of $A+B$ is necessary for
$T$ to be an $\AC(\sigma)$ operator. It would of course be
interesting to know whether it is sufficient.

\begin{cor} Let $T = A+iB$ be a compact $\AC(\sigma)$ operator. Then
$A+B$ is well-bounded.
\end{cor}

\begin{comment}

In general it is an open question whether every $\AC(\sigma)$
operator is actually an $\AC(\sigma(T))$ operator. Answering this
question may depend on determining certain structural properties of
the Banach algebra $\AC(\sigma)$ which would allow general theorems
about non-analytic functional calculi to be applied. We refer the
reader to \cite[Section 1.4]{LM}.

\end{comment}

%%%%%%%%%%%%%%%%%%%%%%%%%%%%%%%%%%%%%%%%%%%%%%%%%%%%%%%%%%%%%%%%%%%%%%%
%
%   Bibliography
%
%%%%%%%%%%%%%%%%%%%%%%%%%%%%%%%%%%%%%%%%%%%%%%%%%%%%%%%%%%%%%%%%%%%%%%%

\end{document}